\documentclass[11pt]{amsart}

\usepackage[english]{babel}
\usepackage[utf8x]{inputenc}
\usepackage{amsmath, amstext, amssymb, amsthm, amscd, graphicx, dsfont, manfnt, mathrsfs, qtree, tipa, shapepar}

\begin{document}

\title{A one-formula proof of the nonvanishing of $L$-functions of real characters at 1 \ \ \ }

\author{Bogdan Veklych}

\begin{abstract}

We present a simple analytic proof that $L$-functions of real non-principal Dirichlet characters are nonzero at 1.

\end{abstract}

\maketitle

It is universally agreed that the most difficult step in proving Dirichlet's theorem on primes in arithmetic progressions lies in showing that the $L$-function, $L(s,\chi)$, attached to a real non-principal Dirichlet character $\chi$ does not vanish at 1. In this note we present a proof of this fact, suitable for the undergraduate classroom. It is similar to those given by Davenport [1] and Serre  [2], but it avoids the use of infinite products and the explicit use of Landau's theorem on Dirichlet series. It does use basic complex analysis and the idea behind Landau's theorem. There are proofs that bypass complex analysis, which are, however, more complicated than the proof below. Of these, a proof by Monsky [3] is perhaps the simplest. Another example is a proof by Apostol [4].  \\

We shall use two standard facts, that the Riemann zeta function $\zeta(s)$ is holomorphic in the half-plane $\operatorname{Re} s>0$ away from $s=1$ where it has a simple pole, and that the $L$-function $L(s, \chi)$ of a non-principal character $\chi$ is holomorphic in all of $\operatorname{Re} s>0$. These facts are easily proved by summation by parts (see, for instance, [1]). We shall assume that there is a real non-principal character $\chi$ for which $L(1,\chi)=0$ and derive a contradiction by considering $F(s)=\zeta(s)L(s, \chi)$. By the cited facts, $F(s)$ is holomorphic in $\operatorname{Re} s >0$ except possibly at $s=1$, but since the zero of $L(s,\chi)$ at 1 cancels the pole of $\zeta(s)$, $F(s)$ is holomorphic everywhere in $\operatorname{Re} s >0$. All we shall need, however, is that $F(s)$ is holomorphic in a disc centered at 2 (or any number greater than 1) that contains $\frac{1}{2}$, or equivalently that $F(2-s)$ is holomorphic in a disc centered at 0 that contains $\frac{3}{2}$.\\

Multiplying for $\operatorname{Re} s>1$ the absolutely convergent Dirichlet series of $\zeta(s)=\sum k^{-s}$ and $L(s, \chi)=\sum \chi(l)l^{-s}$, we obtain a Dirichlet series $\sum c_nn^{-s}$ representing $F(s)$ in (at least) $\operatorname{Re} s>1$. Its coefficients $c_n$ are the Dirichlet convolution of 1 and $\chi$, so $c_n$ is the sum of $\chi(l)$ over all $l$ dividing $n$. Thus, by the multiplicativity of $\chi$, $c_n$ is the product of the contributions of the form $1+\chi(p)+\cdots+\chi(p)^a$ of the primes $p$ dividing $n$, where $a$ is the exponent of $p$ in the factorization of $n$. We use the fact that $\chi(p)$ can take only the values 1, $-1$ and 0 to conclude the following:\\

(1) If $a$ is odd, the contribution of $p$ is $1+a$, 0 or 1;\\
\indent (2) If $a$ is even, the contribution of $p$ is $1+a$ or 1.\\

Hence, $c_n \geq 0$, and if $n$ is a square, $c_n \geq 1$. It follows that the sum $\sum c_nn^{-1/2}$ diverges; by $c_n \geq 0$, it is not less than the similar sum over $n=m^2$, which, by $c_{m^2} \geq 1$, is not less than $1+\frac{1}{2}+\frac{1}{3}+\cdots$. We'll combine the divergence of this sum with our previous observation about the holomorphicity of $F(2-s)$ and several basic facts of complex analysis to get a contradiction.\\

For $s \in [0,1)$, we observe that $F(2-s)$ is equal to the following: \begin{align*} \sum_{n=1}^\infty \frac{c_n}{n^{2-s}}&=\sum_{n=1}^\infty \frac{c_n}{n^2}(e^{s \log n}-1)+\sum_{n=1}^\infty \frac{c_n}{n^2}\\ &=\sum_{n=1}^\infty\sum_{k=1}^\infty\frac{c_n}{n^2}\frac{s^k (\log n)^k}{k!}+F(2)\\&=\sum_{k=1}^\infty s^k\sum_{n=1}^\infty\frac{c_n}{n^2}\frac{(\log n)^k}{k!}+F(2)\end{align*} (the change in the order of summation is permitted since in the double sum all terms are nonnegative). The very last expression is a power series in $s$, and we've just shown that it converges to $F(2-s)$ for all $s \in [0,1)$. Hence the radius of convergence of this power series is at least 1, so it defines a function holomorphic in $|s|<1$. It coincides there with $F(2-s)$ by the uniqueness of the analytic continuation. Therefore, this power series is the power series expansion of $F(2-s)$ around 0, as the latter is unique. By our previous observation that $F(2-s)$ is holomorphic in a disc centered at 0 and containing $\frac{3}{2}$, its power series expansion around 0 must converge to it for all $s$ in this disc, in particular, for $s=\frac{3}{2}$. Thus we can set $s=\frac{3}{2}$ in the power series, and then read the displayed formula backwards, as the change of the order of summation is again permitted. We see that $\sum c_nn^{-1/2}=F\left(\frac{1}{2}\right)$, but $\sum c_nn^{-1/2}$ diverges, as we showed above, which is the contradiction.\\

\end{document}